\documentclass{amsart}
\usepackage{amscd}
\newtheorem{theorem}{Theorem}[section]
\newtheorem{lemma}[theorem]{Lemma}
\newtheorem{proposition}[theorem]{Proposition}
\newtheorem{corollary}[theorem]{Corollary}
\theoremstyle{definition}
\newtheorem{definition}[theorem]{Definition}

\theoremstyle{remark}
\newtheorem{remark}[theorem]{Remark}

\numberwithin{equation}{section}

\theoremstyle{plain}

\newtheorem{problem}[theorem]{\bf Problem}

\def\int{\mathop{\roman{int}}}

\def\diam{\text{diam}}

\def\UUU{{\mathcal U}}
\def\VVV{{\mathcal V}}

\def\diam{\mathrm{diam}}
\def\proof{{\bf Proof. }}
\def\endproof{\hfill \qed}

\numberwithin{equation}{section}
\begin{document}

\title[Topological and uniform structures on universal covering spaces]
{Topological and uniform structures on universal covering spaces}

\author{N.~Brodskiy}
\address[N.~Brodskiy]{Department of Mathematics, University of Tennessee, 227 Ayres Hall, 1403 Circle Drive, Knoxville, TN 37996}
\email{brodskiy@math.utk.edu}

\author{J.~Dydak}
\address[J.~Dydak]{Department of Mathematics, University of Tennessee, 227 Ayres Hall, 1403 Circle Drive, Knoxville, TN 37996}
\email{dydak@math.utk.edu}

\author{B.~Labuz}
\address[B.~Labuz]{Department of Mathematics, University of Tennessee, 227 Ayres Hall, 1403 Circle Drive, Knoxville, TN 37996}
\email{labuz@math.utk.edu}

\author{A.~Mitra}
\address[A.~Mitra]{Department of Mathematics, University of Tennessee, 227 Ayres Hall, 1403 Circle Drive, Knoxville, TN 37996}
\email{ajmitra@math.utk.edu}

\keywords{Universal covering maps, uniform structures}

\subjclass[2000]{Primary 55Q52; Secondary 55M10, 54E15}

\begin{abstract}
We discuss various uniform structures and topologies on the universal covering space $\widetilde X$
and on the fundamental group $\pi_1(X,x_0)$. We introduce a canonical uniform structure $CU(X)$ on a topological space $X$ and use it to relate topologies on $\widetilde X$ and uniform structures on $\widetilde{CU(X)}$.

Using our concept of universal Peano space we show connections between the
topology introduced by  Spanier~\cite{Spa}
and a uniform structure of Berestovskii and Plaut~\cite{BP3}.
We give a sufficient and necessary condition for Berestovskii-Plaut structure to be identical with the one generated by the uniform convergence structure on the space of paths in $X$. We also describe when the topology of Spanier is identical with the quotient of the compact-open topology on the space of paths.
\end{abstract}
\maketitle

\tableofcontents

\section{Introduction}

The classical way of introducing topology on the universal covering space $\widetilde X$
(the space of homotopy classes of based paths in $X$) is as the quotient space
of the space of based paths $(X,x_0)^{(I,0)}$ equipped with the compact-open
topology. We denote the space $\widetilde X$ equipped with this topology by $\widetilde X_{top}$. Spanier introduced a different topology on $\widetilde X$ (see the proof of  Theorem 13 on page 82 of~\cite{Spa}). We call it the whisker topology and use notation $\widetilde X_{wh}$. The whisker topology 
is defined by the basis $B([\alpha],U)=\{[\alpha\ast\beta]\}$, where $\alpha$ is a path in $X$ from $x_0$
to $x_1$, $U$ is a neighborhood of $x_1$ in $X$, and $\beta$ is a path in $U$
originating at $x_1$ (we call the path $\beta$ a {\it $U$-whisker}). This topology was used in 1998 by Bogley and Sieradski~\cite{BogSie}, then
in 2006 by Fisher and Zastrow~\cite{FisZas}. A new topology on $\widetilde X$ called the lasso topology was introduced by the authors~\cite{BDLM3} to characterize the unique path lifting property.

The fundamental group $\pi_1(X,x_0)$ can be identified with the subset of classes of loops of the set $\widetilde X$, thus any topology on  $\widetilde X$ induces a topology on the fundamental group. There have been numerous attempts to introduce topology on the fundamental group~\cite{Hur, Dug, Mel, Bra}. In this paper, when we consider a topology on $\widetilde X$, we mainly address the question of whether the induced topology on the fundamental group makes it a topological group.

There was independent research going on in the realm of uniform structures,
where  Berestovskii and Plaut
\cite{BP3} introduced a new uniform structure on the space of chains
in a uniform space. That structure can be easily generalized to $\widetilde X$,
so we call it the Berestovskii-Plaut uniform structure on $\widetilde X$.
In \cite{BDLM1} the authors used ideas of \cite{BP3}
and Krasinkiewicz-Minc \cite{KraMin} to introduce
generalized uniform paths on a uniform space and to define a uniform
structure $GP(X,x_0)$ on the space of generalized paths originating at $x_0$.
Generalized uniform paths lead naturally to a class of spaces called uniformly
joinable and a theory of generalized uniform covering maps
was developed in \cite{BDLM1} for those spaces. See \cite{BP3} for
a predecessor of that theory.

Uniformly joinable spaces are connected to shape theory
(see \cite{DydSeg} or \cite{MarSeg}) in the sense that pointed $1$-movable continua are always uniformly joinable (see \cite{BDLM1}).
Recently, Fisher and Zastrow \cite{FisZas} found another connection
between shape theory and covering maps. Namely, they showed
the projection $\widetilde X_{wh}\to X$ is a generalized universal covering (in their sense)
if the natural homomorphism $\pi_1(X)\to \check\pi_1(X)$ is a monomorphism
and $X$ is paracompact.
The latter class of spaces includes subsets of surfaces and it is known subcontinua
of surfaces are pointed $1$-movable. This raises the question if one can derive
results of \cite{FisZas} from those in \cite{BDLM1}.
Another idea is to define a new class of covering maps (called Peano covering maps)
for Peano spaces as those spaces form a topological analog of uniformly joinable spaces.

The functor of universal locally-path-connected space of $X$, that we introduce in Section~\ref{SubSection:universalInUniform}
and call the Peano space of $X$, is very useful
as it shows how to derive theorems for arbitrary spaces from analogous
results for locally-path-connected spaces. For example, many results
of \cite{FisZas} follow formally from their special cases involving locally path-connected
spaces only (see~\cite{BDLM3}). A good example is the main result of \cite{FisZas} stating that
the projection $\widetilde X\to X$ is a generalized universal covering if
$X$ is path-connected and the fundemental group of $(X,x_0)$ injects
into its \v Cech fundamental group via the natural homomorphism.

\par It is of interest to investigate relationships between various structures on $\widetilde X$.
One of our first results in that direction is the fact the whisker topology
is philosophically very similar to the quotient compact-open topology.
The difference is that, prior to applying the quotient topology, one takes
the universal locally-path-connected space of $(X,x_0)^{(I,0)}$ (see~\ref{BSTopologyIsQuotientOfP}).
We describe when the whisker topology and the quotient compact-open topology are identical: when the space $X$ is a small loop transfer space. The property of being a strong small loop transfer space characterizes the coinsidence of the whisker topology and the lasso topology.

If $X$ is a uniform space, the set $\widetilde X$ can be equipped with a uniform structure. For basic facts on uniform spaces we refer to \cite{Isb} or \cite{Jam}. We consider three uniform structures on $\widetilde X$: the James uniformity (analogous to the quotient compact-open topology), the Berestovskii-Plaut uniformity (analogous to the whisker topology) and the lasso uniformity. For any connected and uniformly locally path-connected uniform space $X$ the James uniformity on $\widetilde X$ is identical with the lasso uniformity (see~\ref{X_l=X_J}). We also describe when the James uniformity and the Berestovskii-Plaut uniformity are identical: when the space $X$ is a uniform small loop transfer space.

Given a topological space $X$, we show in Section~\ref{SECTION uniform-topology} how to introduce a canonical uniform structure on $X$ that induces the original topology in case $X$ is completely regular. Then we relate various properties (such as local connectivity) of the topological space $X$ to the corresponding uniform properties of the canonical uniform space $CU(X)$. This construction provides a relation between Sections~\ref{SECTION Top-structures-paths} and~\ref{SECTION-U-Structures}: for a topological space $X$ a topological structure on $\widetilde X$ can be obtained as topology induced by a uniform structure on $CU(X)$. For example, the lasso topology on $\widetilde X$ is identical with the topology induced by the lasso uniformity on $\widetilde{CU(X)}$.


\section{Uniform vs topological structures}\label{SECTION uniform-topology}

We will discuss exclusively symmetric subsets $E$ of $X\times X$
(that means $(x,y)\in E$ implies $(y,x)\in E$) and the natural notation here (see \cite{Pla})
is to use $f(E)$ for the set of pairs $(f(x),f(y))$, where $f\colon X\to Y$ is a function.
Similarly, $f^{-1}(E)$ is the set of pairs $(x,y)$ so that $(f(x),f(y))\in E$ if $f\colon X\to Y$
and $E\subset Y\times Y$.
\par The {\bf ball $B(x,E)$ at $x$ of radius $E$} is the set of all $y\in X$
satisfying $(x,y)\in E$. A subset $S$ of $X$ is {\bf $E$-bounded}
if it is contained in $B(x,E)$ for some $x\in X$.
\par
A {\bf uniform structure} on $X$ is a family $\mathcal{E}$ of symmetric subsets $E$
of $X\times X$ (called {\bf entourages}) that contain the diagonal
of $X\times X$, form a filter (that means $E_1\cap E_2\in \mathcal{E}$
if $E_1,E_2\in \mathcal{E}$ and $F_1\in \mathcal{E}$ if $F_2\in \mathcal{E}$
and $F_2\subset F_1$),
and every $G_1\in \mathcal{E}$ admits $G\in \mathcal{E}$ so that
$G^2\subset G_1$ ($G^2$ consists of pairs $(x,z)\in X\times X$
so that there is $y\in X$ satisfying $(x,y)\in G$ and $(y,z)\in G$).
A {\bf base} $\mathcal{F}$ of a uniform structure $\mathcal{E}$
is a subfamily $\mathcal{F}$ of $\mathcal{E}$ so that for every entourage $E$ there is a subset
$F\in \mathcal{F}$ of $E$.

Any filter $\mathcal{F}$ of subsets of $X\times X$ that contain the diagonal
serves as a base of the uniform structure $\mathcal{U}$
defined as all supersets of $\mathcal{F}$ if and only if every $G_1\in \mathcal{F}$ admits $G\in \mathcal{F}$ so that
$G^2\subset G_1$. That is the basic case of generating a uniform
structure.

\par Given a decomposition of a uniform space $X$ the most pressing issue is
if it induces a natural uniform structure on the decomposition space. James \cite[2.13 on p.24]{Jam}
has a concept of weakly compatible relation to address that issue. For the purpose of this paper
we need a different approach.
\begin{definition}\label{GeneratingUCSTructureDef}
Suppose $f\colon X\to Y$ is a surjective function from a uniform space
$X$. {\bf $f$ generates a uniform structure on $Y$} if the family
$f(E)$, $E$ an entourage of $X$, is a base of a uniform structure on $Y$
(that particular uniform structure on $Y$ is said to be {\bf generated by $f$}).
Equivalently, for each entourage $E$ of $X$ there is an entourage $F$ of $X$
such that $f(F)^2\subset f(E)$.
\end{definition}

Notice $f\colon X\to Y$ is uniformly continuous if both $X$ and $Y$ are uniform spaces
and the uniform structure on $Y$ is generated by $f$. Indeed $E\subset f^{-1}(f(E))$
for any entourage $E$ of $X$.

Let us demonstrate how to transfer concepts from topology
to the uniform category:
\begin{definition}\label{UniformPathConnDef}
 A uniform space $X$ is {\bf uniformly locally path-connected}
if for every entourage $E$ of $X$ there is an entourage $F$ of $X$
such that any pair $(x,y)\in F$ can be connected by an $E$-bounded path in $X$.
\par A uniform space $X$ is {\bf uniformly semi-locally simply-connected}
if there is an entourage $F$ of $X$
such that any $F$-bounded loop in $X$ is null-homotopic in $X$.
\end{definition}

\begin{remark}
Notice our definition
of uniform local path-connectedness is much simpler than
\cite[Definition 8.12 on p.119]{Jam} and we are unsure if the definition
\cite[Definition 8.13 on p.119]{Jam} of uniform semi-local simple connectedness
is correct as it involves existence of a base of entourages rather than just one entourage.
\end{remark}

It is well-known (see \cite{Isb} or \cite{Jam})
 there is a forgetful functor from the uniform category
to the topological category. It assigns to a uniform space $X$
the topology defined as follows: $U\subset X$ is open if and only if
for each $x\in U$ there is an entourage $E_x$ of $X$ satisfying
$B(x,E_x)\subset U$. Notice $f\colon X\to Y$ is continuous if
it is uniformly continuous.

Let us describe a functor from the topological category to the uniform
category that is the inverse of the forgetful functor on the class
of completely regular spaces.

\begin{definition}\label{CUSDefinition}
Given a topological space $X$ its {\bf canonical uniform structure}
$CU(X)$ is defined as follows: $E\subset X\times X$ is an entourage
if and only if there is a partition of unity $f=\{f_s\}_{s\in S}$ on $X$
so that $E$ contains 
$$E_f:=\{(x,y)\in X\times X | \text{ there is } s\in S\text{ satisfying }
f_s(x)\cdot f_s(y) > 0\}.$$
\end{definition} 
It is indeed a uniform structure:

\begin{lemma}\label{PULemma}
If $g=\{f'_T\}_{T\subset S}$ is the derivative of $f=\{f_s\}_{s\in S}$,
then $E_g^2\subset E_f$.
\end{lemma}
\proof Recall (see \cite{Dyd}) $f'_T$ is defined as $|T|\cdot \max(0,g_T)$, where
\par\noindent 
$g_T=\min\{f_t\mid t\in T\}-\sup\{f_t\mid t\in S\setminus T\}$ and $T$ is a finite subset of $S$. Suppose $(x,y)\in E_g$ and $(y,z)\in E_g$.
There exist finite subsets $T$ and $Z$ of $S$ such that $f'_T(x)\cdot f'_T(y) > 0$
and $f'_Z(y)\cdot f'_Z(z) > 0$. Hence $T\cap Z\ne\emptyset$
and $f_s(x)\cdot f_s(z) > 0$ for any $s\in T\cap Z$. That means $(x,y)\in E_f$.
\endproof

Notice any continuous function $f\colon X\to Y$ becomes
uniformly continuous when considered as a function from
$CU(X)$ to $CU(Y)$. Indeed, given a partition of unity
$\{f_s\}_{s\in S}$ on $X$, it induces the partition of unity
$\{f\circ f_s\}_{s\in S}$ on $X$.

\begin{proposition}\label{TopologyInducedIsTheSame}
Suppose $X$ is a topological Hausdorff space.
The topology induced on $X$ by $CU(X)$ coincides with that of $X$ if and only if $X$
is completely regular.
\end{proposition}
\proof Notice the topology induced by $CU(X)$ is always coarser than
the original one. Indeed, $B(x,E_f)$ contains the union of
all $f_s^{-1}(0,1]$ so that $f_s(x) > 0$, so every open set in $CU(X)$
belongs to the original topology of $X$.
\par
If $X$ is completely regular and $U$ is open, then for any $x\in X$ there is
a continuous function $g\colon X\to [0,1]$ so that $g(x)=1$ and $g(X\setminus U)\subset \{0\}$. In that case $B(x,E_f)\subset U$, where $f=\{g,1-g\}$. That means $U$
is open in the topology induced by $CU(X)$.
\endproof

\begin{proposition}\label{TopologicalVsUniformConcepts}
Suppose $X$ is a paracompact space.
\begin{itemize}
\item[a.] $X$ is locally connected if and only if $CU(X)$ is uniformly
locally connected.
\item[b.] $X$ is locally path-connected if and only if $CU(X)$ is uniformly
locally path-connected.
\item[c.] If $X$ is locally path-connected, then it is semi-locally simply-connected if and only if $CU(X)$ is uniformly
semi-locally simply-connected.
\end{itemize}
\end{proposition}
\proof c. Suppose $X$ is locally path-connected.
Choose an open cover $\{V_x\}_{x\in X}$ of $X$ with the property that any loop in $V_x$ based at $x$ is null-homotopic
in $X$. Let $W_x$ be the path component of $V_x$ containing $x$.
Observe that any loop in $W_x$ is null-homotopic in $X$.
Choose a partition of unity $g=\{g_x\}_{x\in X}$  on $X$ so that
$g_x^{-1}(0,1]\subset V_x$ for all $x\in X$. 
Let $f$ be its derivative. 
Suppose $\alpha$ is a loop at $x$ in $B(x,E_f)$. Observe as in \ref{PULemma}
that
there is $y\in X$ so that $B(x,E_f)\subset W_y$. Therefore $\alpha$
is null-homotopic in $X$.
\endproof

\section{Universal spaces in the uniform category}\label{SubSection:universalInUniform}
In this section we generalize concepts from~\cite{BDLM3}
to the uniform category. 

\begin{definition}\label{UPeanoSpacesDef}
A uniform space $X$ is a {\bf uniform Peano space}
if it is uniformly locally path-connected and connected.
\end{definition}

In analogy to the universal covering spaces we introduce the following notion:

\begin{definition}\label{ULPCspacesDef}
Given a uniform space $X$ its {\bf universal Peano space}
$P(X)$ is a uniform Peano space together with a uniformly continuous map
$\pi\colon P(X)\to X$ satisfying the following universality condition:
\par For any map $f\colon Y\to X$ from a uniform Peano space $Y$
there is a unique uniformly continuous lift $g\colon Y\to P(X)$ of $f$ (that means $\pi\circ g=f$).
\end{definition}

Notice $P(X)$ is unique if it exists.

\begin{theorem}\label{ULPCExistsThm}
Every path-connected uniform space $X$ has a universal Peano space. It is the set $X$ equipped
with the uniform structure generated by $\{pc(E)\}_E$, where $pc(E)$ is the set of
pairs $(x,y)$ that can be connected by an $E$-bounded path in $X$, $E$ an entourage of $X$.
\end{theorem}
\proof Notice $(pc(F))^4\subset pc(E)$ if $F^2\subset E$, so $\{pc(E)\}_E$ does indeed generate
a uniform structure on the set $X$ resulting in a new uniform space $P(X)$
so that the identity function $P(X)\to X$ is uniformly continuous.
\par Given a uniformly continuous function $f\colon Y\to X$ from a uniform Peano space $Y$ and given an entourage $E$ of $X$ there is an entourage $F$ of $Y$
contained in $f^{-1}(E)$ with the property that any $F$-close points in $Y$
can be connected by an $f^{-1}(E)$-bounded path. Therefore $F\subset f^{-1}(pc(E))$ which proves $f\colon Y\to P(X)$
is uniformly continuous. It also proves $P(X)$ is path-connected as any path in $X$ induces a path in $P(X)$. Since any $pc(E)$-close points can be connected by $E$-bounded path, $P(X)$ is uniformly locally path-connected.
\endproof

\begin{corollary}\label{MetrizabilityOfP}
If $X$ is a metrizable uniform space, then so is its universal Peano space $P(X)$.
\end{corollary}
\proof Notice $P(X)$ is Hausdorff and has a countable base of entourages. Such spaces are metrizable by 
 \cite[Theorem 13 on p.186]{Kel}.
\endproof

\begin{theorem}\label{UNIversalCechSpace}
For every pointed uniform Peano space $(X,x_0)$ there is a universal object
in the class of uniformly continuous maps $f\colon (Y,y_0)\to (X,x_0)$ such that $Y$ is a uniform Peano space and $f$ induces the trivial homomorphism
$\check \pi_1(f)\colon \check \pi_1(Y,y_0)\to \check \pi_1(X,x_0)$ of the uniform
fundamental groups. 
\end{theorem}
\proof It is shown in \cite{BDLM1} that there is a uniformly continuous
map $\pi_X\colon GP(X,x_0)\to (X,x_0)$ from a uniform space of the trivial uniform fundamental group with the property
that for any uniformly continuous map $f\colon (Y,y_0)\to (X,x_0)$
from a pointed uniform Peano space $(Y,y_0)$ there is a unique lift
$g\colon (Y,y_0)\to GP(X,x_0)$ if and only if $f$ induces
the trivial homomorphism of uniform fundamental groups (the result in \cite{BDLM1}
is for so-called uniformly joinable spaces and it applies here as uniform Peano spaces
are uniformly joinable).
\par Let $(Z,z_0)$ be the universal Peano space of $GP(X,x_0)$
 with the induced
uniformly continuous map $\pi\colon (Z,z_0)\to (X,x_0)$.
Notice any uniformly continuous function $f\colon (Y,y_0)\to (X,x_0)$ 
has a unique lift to $(Z,z_0)$ provided it induces the trivial homomorphism
of the uniform fundamental groups. Since $\pi\colon (Z,z_0)\to (X,x_0)$
factors through $GP(X,x_0)$ it does induce the trivial homomorphism
of the uniform fundamental groups.
\endproof

In order to apply \ref{UNIversalCechSpace} we need to
discuss simple-connectivity of the space $GP(X,x_0)$ of generalized paths
in $X$ originating from $x_0$.

\begin{proposition}\label{SimpleConnectivityOfGP}
Suppose $(X,x_0)$ is a pointed uniform Peano space.
$GP(X,x_0)$ is simply connected if and only if
the natural homomorphism $\pi_1(X,x_0)\to \check \pi_1(X,x_0)$ is a monomorphism.
\end{proposition}
\proof Suppose $\pi_1(X,x_0)\to \check \pi_1(X,x_0)$ is a monomorphism.
Given a loop $\alpha\colon (S^1,1)\to GP(X,x_0)$, its composition
with $\pi_X\colon GP(X,x_0)\to (X,x_0)$ is trivial in $\check\pi_1(X,x_0)$
as $\check\pi_1(\pi_X)=0$. Therefore $\pi_X\circ\alpha$ is trivial in $\pi_1(X,x_0)$
and can be extended over the $2$-disk $D^2$.
That map has a unique lift to $GP(X,x_0)$, hence $\alpha$ is null-homotopic.
\par Suppose $GP(X,x_0$ is simply connected and $\beta\colon (S^1,1)\to (X,x_0)$
becomes trivial in $\check\pi_1(X,x_0)$. In that case $\beta$ lifts to $GP(X,x_0)$
and is null-homotopic.
\endproof

\begin{corollary}\label{PiUniversalSpaceExistsThm}
Suppose $(X,x_0)$ is a pointed uniform space. If the natural homomorphism
$\pi_1(X,x_0)\to \check \pi_1(X,x_0)$ is a monomorphism,
then $(X,x_0)$ has a universal simply-connected Peano space.
\end{corollary}
\proof From the commutativity of
$$
\begin{CD}
\pi_1(P(X,x_0))  @> >>  \pi_1(X,x_0)    \\
@V VV                               @VV V \\
 \check\pi_1(P(X,x_0))   @> >>   \check\pi_1(X,x_0)
\end{CD}
$$
we get the natural homomorphism $ \pi_1(P(X,x_0))\to  \check\pi_1(P(X,x_0))$
is a monomorphism and by~\ref{ULPCExistsThm} we can reduce the general case to that of $(X,x_0)$
being a pointed uniform Peano space.
\par

Consider $SP(X,x_0)= P(GP(X,x_0))$ with the induced map 
\par\noindent $\pi\colon
SP(X,x_0)\to (X,x_0)$. Since $\check\pi_1(SP(X,x_0))\to \check\pi_1(X,x_0)$ is trivial
and $\pi_1(X,x_0)\to \check \pi_1(X,x_0)$ is a monomorphism,
$\pi_1(SP(X,x_0))\to \pi_1(X,x_0)$ is trivial.
\par Given any uniformly continuous $f\colon (Y,y_0)\to (X,x_0)$ from a pointed uniform Peano space
so that $\pi_1(f)=0$ it suffices to show
 $\check \pi_1(f)\colon \check\pi_1(Y,y_0)\to\check\pi_1(P(X,x_0))$ is trivial
 as then  it lifts uniquely to $P(GP(X,x_0))$.
 \par Given an entourage $E$ of $P(X)$ pick an entourage $F$ of $Y$
 such that any two points $(x,y)\in F$ can be connected by
 a path contained in $B(x,f^{-1}(E))$.
 Given an element $\gamma$ of $\check\pi_1(Y,y_0)$ its representative
 in $\mathcal{R}(Y,F)$ leads to a loop $\gamma_E$ in $Y$ at $y_0$
 that represents $\gamma$ in $\mathcal{R}(Y,f^{-1}(E))$.
 As $f(\gamma_E)$ is null-homotopic in $P(X)$,
 $\check\pi_1(f)(\gamma)=1$.
\endproof

\begin{proposition}\label{GPVsBP}
Suppose $(X,x_0)$ is a pointed uniform space. If the natural homomorphism
$\pi_1(X,x_0)\to \check \pi_1(X,x_0)$ is a monomorphism,
then the natural function from $\widetilde X_{BP}$ to $P(GP(P(X,x_0)))$
is a uniform equivalence.
\end{proposition}
\proof Since each path in $(X,x_0)$ is a generalized path in $P(X,x_0)$,
there is a natural function $i\colon \widetilde X_{BP}\to GP(P(X,x_0))$.
It is uniformly continuous and $\widetilde X_{BP}$ is a uniform Peano space,
hence $i\colon \widetilde X_{BP}\to P(GP(P(X,x_0)))$ is uniformly continuous.
Since $P(GP(P(X,x_0)))$ is path-connected, its elements are generalized paths in 
$P(X,x_0)$ that are representable by genuine paths. That means $i$
is surjective.
If $i$ is injective, then any loop in $(X,x_0)$ 
\endproof

\begin{corollary}[Fisher-Zastrow \cite{FisZas}]\label{FisZasTheorem}
Suppose $X$ is a path-connected paracompact space. If the natural homomorphism
$\pi_1(X,x_0)\to \check \pi_1(X,x_0)$ is a monomorphism for some $x_0\in X$,
then the projection $\widetilde X_{wh}\to X$ of $X$ is a generalized universal cover.
\end{corollary}
\proof Consider $X$ equipped with the canonical uniform structure $CU(X)$.
As $\check\pi_1(CU(X),x_0)=\check\pi_1(X,x_0)$,
there is a universal map $f\colon Y\to CU(X)$
among all uniformly continuous functions from uniform Peano spaces
that induce trivial homomorphism of fundamental groups. 
Given a Peano space $Z$ and a map $g\colon (Z,z_0)\to (X,x_0)$
inducing trivial homomorphism of the fundamental groups, we consider
$Z$ equipped with the canonical uniform structure $CU(Z)$.
Notice in the proof of \ref{PiUniversalSpaceExistsThm}
that we only need to use $CU(Z)$ is uniformly joinable to lift $g$
to $GP(X,x_0)$ (i.e., we do not need $CU(Z)$ to be a uniform Peano space
which it certainly is if $Z$ is paracompact by~\ref{TopologicalVsUniformConcepts}). 
Thus the lift exists, it lands in the path component of the trivial generalized path,
and then we can lift it to the Peano space of that component as $Z$ is a Peano space.

Consider $SP(X,x_0)= P(GP(P(X,x_0)))$ with the induced map 
\par\noindent $\pi\colon
SP(X,x_0)\to (X,x_0)$. Since $\check\pi_1(SP(X,x_0))\to \check\pi_1(X,x_0)$ is trivial
and $\pi_1(X,x_0)\to \check \pi_1(X,x_0)$ is a monomorphism,
$\pi_1(SP(X,x_0))\to \pi_1(X,x_0)$ is trivial.
\par Given any uniformly continuous $f\colon Y\to X$ from a uniform Peano space
so that $\pi_1(f)=0$ its lift $f\colon Y\to P(X)$ also satisfies $\pi_1(f)=0$,
so it suffices to show
 $\check \pi_1(f)\colon \check\pi_1(Y,y_0)\to\check\pi_1(P(X,x_0))$ is trivial
 as then  it lifts uniquely to $P(GP(P(X,x_0)))$.
 \par Given an entourage $E$ of $P(X)$ pick an entourage $F$ of $Y$
 such that any two points $(x,y)\in F$ can be connected by
 a path contained in $B(x,f^{-1}(E))$.
 Given an element $\gamma$ of $\check\pi_1(Y,y_0)$ its representative
 in $\mathcal{R}(Y,F)$ leads to a loop $\gamma_E$ in $Y$ at $y_0$
 that represents $\gamma$ in $\mathcal{R}(Y,f^{-1}(E))$.
 As $f(\gamma_E)$ is null-homotopic in $P(X)$,
 $\check\pi_1(f)(\gamma)=1$.
\endproof

\section{Topological structures on $\widetilde X$}\label{SECTION Top-structures-paths}

In this section we discuss topological structures on $\widetilde X$
and $\pi_1(X,x_0)$. Notice that neither of the two topologies on $\pi_1(X,x_0)$ discussed in this section makes $\pi_1(X,x_0)$ a topological group~\cite{BDLM3}.

In case of $X$ being metric $\pi_1(X,x_0)$
is considered as a quotient space $\pi_1^{top}(X,x_0)$ of the space
$(X,x_0)^{(S^1,1)}$ of loops equipped with the uniform convergence
metric (see \cite{Fab}).
We generalize that concept as follows:

\begin{definition}\label{QuotientCODef}
For any pointed topological space $(X,x_0)$ let $\widetilde X_{top}$ be the quotient space
of the function space $(X,x_0)^{(I,0)}$ equipped with the compact-open topology.
$\pi_1^{top}(X,x_0)$ is the quotient space of $(X,x_0)^{(S^1,1)}$ equipped with the compact-open topology.
\end{definition}

Spanier~\cite{Spa} introduced a new topology on $\pi_1(X,x_0)$ that was used by Bogley and Sieradski~\cite{BogSie}. This topology was generalized later by Fisher and Zastrow~\cite{FisZas} to $\widetilde X$ and used by the authors~\cite{BDLM3}. 

\begin{definition}\label{BSTopDef}
For any pointed topological space $(X,x_0)$ the {\it whisker} topology on the set $\widetilde X$
is defined by the basis $B([\alpha],U)=\{[\alpha\ast\beta]\}$, where $\alpha$ is a path in $X$ from $x_0$
to $x_1$, $U$ is a neighborhood of $x_1$ in $X$, and $\beta$ is a path in $U$
originating at $x_1$ (we call the path $\beta$ a {\it $U$-whisker}). We denote $\widetilde X$ with the whisker topology by $\widetilde X_{wh}$.
By $\pi_1^{wh}(X,x_0)$ we mean $\pi_1(X,x_0)$ with the subspace
topology inherited from $\widetilde X_{wh}$.
\end{definition}

Given a neighborhood $U$ of $\alpha(1)$ in $X$
and a neighborhood $W$ of $x_0$ in $X$
define $B([\alpha],W,U)$ as elements of $\widetilde X$
of the form $\beta\ast\alpha\ast\gamma$, where $\beta$ is a loop at $x_0$
in $W$ and $\gamma$ is a path in $U$ originating at $\alpha(1)$.

Suppose $\alpha$ and $\beta$ are two paths in $X$
so that $\alpha(1)=\beta(0)$. Given a neighborhood $U_0\times U_1\times U_2$ of 
$(\alpha(0),\alpha(1),\beta(1))$ in $X^3$,
we define $B([\alpha\ast\beta],U_0,U_1,U_2)$
as the set of homotopy classes (rel. end-points) of paths
of the form $\omega\ast\alpha\ast\gamma\ast\beta\ast\lambda$,
where $\omega$ is a loop in $U_0$ at $\alpha(0)$,
$\gamma$ is a loop in $U_1$ at $\alpha(1)$, and $\lambda$
is a path in $U_2$ originating at $\beta(1)$.

\begin{lemma}\label{BasicLemma}
Let $\pi\colon (X,x_0)^{(I,0)}\to \widetilde X$ be the projection.
If $V$ is a neighborhood of the concatenation
$\alpha\ast\beta$ of two paths in $ (X,x_0)^{(I,0)}$,
then there is a neighborhood $U_0\times U_1\times U_2$ of 
$(x_0,\alpha(1),\beta(1))$ in $X^3$
such that $B([\alpha\ast\beta],U_0,U_1,U_2)\subset \pi(V)$.
\end{lemma}
\proof There exist a finite sequence
of pairs $\{(C_s,V_s)\}_{s\in S}$ such that $\alpha(C_s)\subset V_s$ for all
$s\in S$, where $C_s$ are compact subsets of $I$, $V_s$ are open subsets of $X$,
and any path $\beta$ in $X$ satisfying $\beta(C_s)\subset V_s$ for all $s\in S$
belongs to $V$. Given $t\in I$, Let $U_{2t}$ be the intersection
of $\{V_s\}_{s\in T}$, where $T$ consists of $s\in S$ so that $(\alpha\ast\beta)(t)\in C_s$
(if $T$ is empty, we put $U_{2t}=X$). 

Suppose $\omega$ is a loop in $U_0$ at $\alpha(0)$,
$\gamma$ is a loop in $U_1$ at $\alpha(1)$, and $\lambda$
is a path in $U_2$ originating at $\beta(1)$.

Given $s< t$ let $u_{s,t}\colon [s,t]\to [0,1]$ be the increasing linear homeomorphism. Given $t < \frac{1}{4}$, define $\mu_t\colon [0,1]\to X$ as follows:

\begin{equation*}
\mu_t(x)= 
\begin{cases} 
\omega(x) & \text{if $x\leq t$,}\\
\alpha(x) & \text{if $t\leq x\leq 1/2-t$,}\\
\gamma(x) & \text{if $1/2-t\leq x\leq 1/2+t$,}\\
\beta(x) & \text{if $1/2+t\leq x\leq 1-t$,}\\
\lambda(x) & \text{if $1-t\leq x\leq 1$.}\\
\end{cases}
\end{equation*}

Notice $\mu_t(C_s)\subset V_s$ for $t$ sufficiently close to $0$.
For such $t$ we have $\mu_t\in V$, hence $[\omega\ast\alpha\ast\gamma\ast\beta\ast\lambda]=[\mu_t]\in \pi(V)$.
That proves $B([\alpha\ast\beta],U_0,U_1,U_2)\subset \pi(V)$.
\endproof

\begin{theorem}\label{BSTopologyIsQuotientOfP}
For any pointed topological space $(X,x_0)$ its universal covering space
$\widetilde X$ equipped with the whisker topology is identical
with the quotient space of the universal Peano space of $(X,x_0)^{(I,0)}$.
\end{theorem}
\proof Let $\pi\colon (X,x_0)^{(I,0)}\to \widetilde X$ be the projection. Suppose $\beta\in\pi^{-1}(B([\alpha],U))$ for some paths $\alpha,\beta\in (X,x_0)^{(I,0)}$
such that $U$ is an open neighborhood of the end-point of $\alpha$.
That means $\beta\sim \alpha\ast\gamma$ for some path $\gamma$ in $U$
starting from $\alpha(1)$. Consider $V=\{\omega\in (X,x_0)^{(I,0)} \mid \omega(1)\subset U\}$
and assume $\lambda$ belongs to the path-component of $V$ containing $\beta$.
That implies $\lambda\sim \beta\ast \mu$, where $\mu$ is a path in $U$.
Thus $\lambda\sim \alpha\ast(\gamma\ast\mu)$ and $[\lambda]\in B([\alpha],U)$.
That proves $P((X,x_0)^{(I,0)})\to \widetilde X_{wh}$ is continuous.
\par Suppose $\pi^{-1}(A)$ is open in $P((X,x_0)^{(I,0)})$ for some
subset $A$ of $\widetilde X$ and $\alpha\in \pi^{-1}(A)$. There is an open neighborhood $V$ 
of $\alpha$ whose path component containing $\alpha$ is a subset of $\pi^{-1}(A)$.
By \ref{BasicLemma} there is a neighborhood $U$ of $\alpha(1)$ in $X$
 so that $B([\alpha],U)\subset A$. As $B([\alpha],U)$ is path-connected,
 $A$ is open.
 \endproof

\begin{problem}\label{BSTopologyInTwoWays}
Let $(X,x_0)$ be a pointed topological space.
Is $\widetilde X_{wh}$ the universal Peano space of $\widetilde X_{top}$?
Equivalently, is the identity function $P(\widetilde X_{top})\to \widetilde X_{wh}$
continuous?
\end{problem}

It is of interest to determine when the compact-open topology coincides with the whisker
topology on $\widetilde X$. It turns out (see \ref{PiContinuousImpliesSLT}) it is related
to continuity
of the induced homomorphism
$h_\alpha\colon \pi_1^{wh}(X,x_1)\to \pi_1^{wh}(X,x_0)$
($h_\alpha([\beta]):=\alpha^{-1}\ast\beta\ast\alpha$), where $\alpha$ is a path
from $x_0$ to any $x_1\in X$. Therefore it makes sense to give a necessary and sufficient condition for continuity of $h_\alpha$.

\begin{lemma}\label{ContinuityOfhalpha}
$h_\alpha\colon \pi_1^{wh}(X,x_1)\to \pi_1^{wh}(X,x_0)$
is continuous if and only if for any neighborhood $U$ of $x_0$
there is a neighborhood $V$ of $x_1$ such that for any loop $\beta$
in $V$ based at $x_1$, the loop
$\alpha\ast\beta\ast\alpha^{-1}$
is homotopic rel.~$x_0$ to a loop contained in $U$.
\end{lemma}
\proof Suppose $h_\alpha$ is continuous and $U$ is a neighborhood of $x_0$.
There is a neighborhood $V$ of $x_1$ such that $h_\alpha(B(1,V))\subset
B(1,U)$ (here $1$ is used for the trivial loops at $x_0$ and $x_1$).
Given any loop $\gamma$ in $V$ at $x_1$ one has $[\gamma]\in B(1,V)$,
so $h_\alpha([\gamma])\in B(1,U)$.
\par
Assume $h_\alpha(B(1,V))\subset
B(1,U)$.
Given any loop $\gamma$ in $V$ at $x_1$ and given any loop $\omega$ at $x_1$
one has $h_\alpha([\omega\ast\gamma])=
h_\alpha([\omega])\ast h_\alpha([\gamma])\in B(h_\alpha([\omega]),U)$,
so $h_\alpha$ is continuous.
\endproof

\begin{definition}\label{TopologicalSLTDef}
A topological space $X$ is a {\bf small loop transfer space} (SLT-space for short)
if for every path $\alpha$ in $X$ and every neighborhood $U$ of $x_0=\alpha(0)$
there is a neighborhood $V$ of $x_1=\alpha(1)$
such that given a loop $\beta\colon (S^1,1)\to (V,x_1)$
there is a loop $\gamma\colon (S^1,1)\to (U,x_0)$ that is homotopic
to $\alpha\ast\beta\ast\alpha^{-1}$ rel.~$x_0$.
\end{definition}

\begin{problem} \label{AreSLTsSimplyConnected}
Is there a Peano SLT-space that is not semi-locally simply
connected?
\end{problem}

Let us explain SLT-spaces in terms of homomorphism
$h_\alpha\colon \pi_1^{wh}(X,x_1)\to \pi_1^{wh}(X,x_0)$ induced by a path
$\alpha$ in $X$ from $x_0$ to $x_1$ ($h_\alpha([\beta])=\alpha\ast\beta\ast\alpha^{-1}$):

\begin{corollary}\label{SLTMeaning}
$X$ is an SLT-space if and only if every path $\alpha$ in $X$ from $x_0$ to $x_1$
induces a homeomorphism $h_\alpha\colon \pi_1^{wh}(X,x_1)\to \pi_1^{wh}(X,x_0)$
\end{corollary}

\begin{proposition}\label{HEIsNotSLT}
The Hawaiian Earring is not an SLT-space.
\end{proposition}
\proof The Hawaiian Earring $HE$ is the union of infinitely
many circles $C_n$ tangent to each other at the point $0$
so that $\diam(C_n)=\frac{1}{n}$.
Let $g_n$ be the loop defined by circle $C_n$.
If $HE$ is an SLT-space, then
there is $m > 0$ such that  $g_1^{-1}\ast g_m\ast g_1$ is 
 homotopic rel. end-points to a loop $\alpha$
of diameter less that $1$.
Consider the retraction $r\colon HE\to C_1\cup C_m$
sending $C_i$ homeomorphically onto $C_m$ if $1< i\leq m$
and sending $C_i$ to $0$ if $i > m$.
Applying $r$ to $g=g_1^{-1}\ast g_m\ast g_1$ shows $g$ is homotopic
rel. end-points to a power of $g_m$, a contradiction.
\endproof

\begin{theorem}\label{PiContinuousImpliesSLT}
If $\widetilde X_{top}=\widetilde X_{wh}$ for all $x_0\in X$,
then $X$ is an SLT-space.
\end{theorem}
\proof Consider the projections $\pi\colon (X,x_0)^{(I,0)}\to \widetilde X_{wh}$
from 
\par\noindent $(X,x_0)^{(I,0)}$ equipped with the compact-open topology.
By \ref{BSTopologyIsQuotientOfP} the projections $\pi$ are continuous
for all $x_0\in X$.

Suppose $\alpha$ is a path in $X$ from $x_0$ to $x_1$ and $U$ is a neighborhood
of $x_1$. Since $\pi^{-1}(B([\alpha],U))$ is open and contains $\alpha$,
\ref{BasicLemma} says there are neighborhoods $V$ of $x_0$ in $X$
and $W$ of $x_1$ in $X$ such that $B([\alpha],V,W)\subset B([\alpha],U)$.
Suppose $\beta$ is a loop in $V$ at $x_0$.
There is a loop $\gamma$ in $U$
at $x_1$ so that $\beta\ast\alpha$ is homotopic to $\alpha\ast\gamma$
rel. end-points. That is equivalent to $X$ being an SLT-space.
\endproof

\begin{theorem}\label{SLTImpliesPiContinuous}
If $X$ is a Peano SLT-space, then
$\pi_1^{wh}(X,x_0)=\pi_1^{top}(X,x_0)$
and $\widetilde X_{top}=\widetilde X_{wh}$.
\end{theorem}
\proof Suppose $\alpha$ is a path from $x_0$ to $x_1$ and $U$ is a neighborhood
of $x_1$ in $X$. For each $t\in [0,1]$ let $\alpha_t$ be the path
from $x_t=\alpha(t)$ to $x_1$ determined by $\alpha$.
For such $t$ choose a path-connected neighborhood $V_t$ of $x_t$
with the property that for any loop $\beta$ in $V_t$ at $x_t$
there is a loop $\gamma$ in $U$ at $x_1$ so that $\beta\ast\alpha_t$ is homotopic rel. end-points
to $\alpha_t\ast\gamma$.
\par For each $t\in I$ choose a closed subinterval $I_t$ of $I$ containing $t$ in its interior
(rel.~$I$) so that $\alpha(I_t)\subset V_t$. Choose finitely many
of them that cover $I$ so that no proper subfamily covers $I$.
Let $S$ be the set of such chosen points $t$.
If $I_s\cap I_t\ne\emptyset$, let $V_{s,t}$ be the path component of $V_s\cap V_t$
containing $\alpha(I_s\cap I_t)$.
Otherwise, $V_{s,t}=\emptyset$.
Let $W$ be the set of all paths $\beta$ in $X$ originating at $x_0$ such that
$\beta(I_s\cap I_t)\subset V_{s,t}$ for all $s,t\in S$.
It suffices to show $\pi(W)\subset B([\alpha],U)$.
\par Given $\beta$ in $W$ pick points $x_{s,t}\in I_s\cap I_t$ if $s\ne t$
and $I_s\cap I_t\ne\emptyset$. Then choose paths $\gamma_{s,t}$ in $X$
from $\alpha(x_{s,t})$ to $\beta(x_{s,t})$.
Arrange points $x_{s,t}$ in an increasing sequence $y_i$, $1\leq i\leq n$,
put $y_0=0$, $y_{n+1}=1$,
and create loops $\lambda_i$ in $X$ at $\alpha(y_i)$ as follows:
travel along $\alpha$ from $\alpha(y_i)$ to $\alpha(y_{i+1})$, then along $\gamma_{i+1}$,
then reverse $\beta$ from $\beta(y_{i+1})$ to $\beta(y_i)$, finally reverse $\gamma_i$.
Notice $\alpha^{-1}\ast\beta$ is homotopic to
$\prod \alpha_i^{-1}\ast \lambda_i\ast\alpha_i$, where $\alpha_i$ is the path
determined by $\alpha$ from $\alpha(y_i)$ to $x_1$.
Each of $ \alpha_i^{-1}\ast \lambda_i\ast\alpha_i$ is homotopic rel.~$x_1$ 
to a loop in $U$, so we are done. 
\endproof

In connection to \ref{AreSLTsSimplyConnected} let us introduce a subclass of SLT-spaces:

\begin{definition}\label{SmallLoopSpacesDef}
A topological space $X$ is called a {\bf small loops space} (or an {\bf SL-space} for short)
if for every loop $\alpha$ in $X$ at $x$ and for any neighborhood $U$ of $x$
there is a loop in $U$ homotopic rel. end-points to $\alpha$.
\end{definition}

SL-spaces are completely opposite to {\bf homotopically Hausdorff spaces}
that are defined by the condition that for any non-null-homotopic loop
$\alpha$ at any $x\in X$ there is a neighborhood $U_\alpha$ of $x$
such that $[\alpha]$ does not belong to the image of
$\pi_1(U_\alpha,x)\to \pi_1(X,x)$
in 
(see \cite{FisZas}). 

The following problem was solved by \v{Z}iga Virk~\cite{V1, V2} after the preprint of this paper was posted.

\begin{problem}\label{AreSLSimplyConnected}
Is there an SL-space $X$ that is not simply-connected?
\end{problem}

We introduced the lasso topology on the set $\widetilde X$ in~\cite{BDLM3}. 

\begin{definition}
Let $\UUU$ be an open cover of a topological space $X$ and $x$ be a point in $X$. A path $l$ is called {\it $\UUU$-lasso based at the point $x$} if $l$ is equal to a finite concatenation of loops $\alpha_n\ast\gamma_n\ast\alpha_n^{-1}$,
where $\gamma_n$ is a loop in some $U\in\mathcal{U}$
and $\alpha_n$ is a path from $x$ to $\gamma_n(0)$.
\end{definition}

\begin{remark}\label{ConjugateLassoIsLasso}
Any conjugate of a $\UUU$-lasso is a $\UUU$-lasso. Indeed, $$\beta*\left(\prod_{i=1}^n \alpha_i\ast\gamma_i\ast\alpha_i^{-1}\right)*\beta^{-1}=\prod_{i=1}^n \left(\beta*\alpha_i\ast\gamma_i\ast\alpha_i^{-1}*\beta^{-1}\right)=\prod_{i=1}^n \left(\beta*\alpha_i\right)\ast\gamma_i\ast\left(\beta*\alpha_i\right)^{-1}$$
\end{remark}

\begin{definition}
For any topological space $X$ the {\it lasso} topology on the set $\widetilde X$ is defined by the basis 
$B([g],\UUU,W)$, where $W$ is a neighborhood of the endpoint $g(1)$ and $\UUU$ is an open cover of $X$. A homotopy class $[h]\in \widetilde X$ belongs to $B([g],\UUU,W)$ if and only if this class has a representative of the form $l * g * \beta$, where $l$ is a $\UUU$-lasso based at $g(0)$ and $\beta $ is a $W$-whisker of $g$.
We denote by $\widetilde X_l$ the set $\widetilde X$ equipped with the lasso topology.
\end{definition}

\begin{definition}\label{TopologicalStrongSLTDef}
A topological space $X$ is a {\bf strong small loop transfer space} (strong SLT-space for short)
if for each point $x_1$, each neighborhood $U$ of $x_1$, and each point $x_2$ there is a neighborhood $V$ of $x_2$ such that for every path $\alpha$ in $X$ from $x_1$ to $x_2$ given any loop $\beta\colon (S^1,1)\to (V,x_2)$
there is a loop $\gamma\colon (S^1,1)\to (U,x_1)$ that is homotopic
to $\alpha\ast\beta\ast\alpha^{-1}$ rel.~$x_1$.
\end{definition}

\begin{proposition}\label{lasso=whiskerIFFstrongSLT}
A path-connected topological space $X$ is a strong SLT-space if and only if $\widetilde X_{l}=\widetilde X_{wh}$ for all $x_0\in X$.
\end{proposition}
\proof
Suppose $X$ is a strong SLT-space. Given any point $x_0\in X$ and a basic neighborhood $B([\omega],U)$ in $\widetilde X_{wh}$, we need to find a basic neighborhood $B([\omega],\VVV,W)$ of $\widetilde X_{l}$ inside $B([\omega],U)$. Any point $x\in X$ has a neighborhood $V$ defined by the strong SLT-property applied to the point $x_1=\omega(1)$ and its neighborhood $U$. Let $\VVV$ be a cover of $X$ by such neighborhoods $V$. Put $W=U$ and consider any element $[l*\omega *\beta] \in B([\omega],\VVV,W)$. Suppose $l$ is equal to a finite concatenation of loops $l=\prod_{i=1}^n \alpha_i\ast\gamma_i\ast\alpha_i^{-1}$, where $\gamma_i$ is a loop in some $V\in\VVV$
and $\alpha_i$ is a path from $x_0$ to $\gamma_i(0)$. Then $[l*\omega *\beta]=[\omega*\left(\prod_{i=1}^n \omega^{-1}*\alpha_i\ast\gamma_i\ast\alpha_i^{-1}*\omega\right)*\beta]$. By the strong SLT property each loop $\omega^{-1}*\alpha_i\ast\gamma_i\ast\alpha_i^{-1}*\omega$ is homotopic rel.~$x_1$ to a loop inside $U$, thus $[l*\omega *\beta]$ belongs to $B([\omega],U)$.

Suppose now that $\widetilde X_{l}=\widetilde X_{wh}$ for any choice of the basepoint.
Let $x_1\in X$ be a point and $U$ be its neighborhood. Denote by $[x_1]$ the trivial path based at $x_1$. Consider the spaces $\widetilde X_{l}=\widetilde X_{wh}$ based at $x_1$. The neighborhood $B([x_1],U)$ in $\widetilde X_{wh}$ contains a lasso neighborhood $B([x_1],\VVV,W)$, where $W\subset U$. Let $x_2$ be any point in $X$ and $V$ be any element of the cover $\VVV$ containing $x_2$. Then for any path $\alpha$ in $X$ from $x_1$ to $x_2$ and any loop $\beta\colon (S^1,1)\to (V,x_2)$ the class $[\alpha\ast\beta\ast\alpha^{-1}]\in B([x_1],\VVV,W)\subset B([x_1],U)$ has a representative $\gamma\colon (S^1,1)\to (U,x_1)$.
\endproof

Since we introduce topologies on the fundamental group, it is natural to ask if these topologies make the fundamental group a topological group (i.e. if the operations of multiplication and taking inverse are continuous in the given topology). The fundamental group with the lasso topology is a topological group~\cite{BDLM3}. In general, neither the quotient of the compact-open topology nor the whisker topology make the fundamental group a topological group~\cite{Fab3, BDLM3}. It is interesting to note that these topologies fail to make $\pi_1(X,x_0)$ a topological group for different reasons. 
The inversion is always continuous in the group $\pi^{top}_1(X,x_0)$ (which makes it a quasi topological group~\cite{AT}).
On the other hand, the continuity of the inversion in the group $\pi^{wh}_1(X,x_0)$ would make it a topological group:

\begin{proposition}\label{inverseWHcontinuous}
Let $(X,x_0)$ be a pointed topological space.
If the operation of taking inverse in $\pi^{wh}_1(X,x_0)$ is continuous, then $\pi^{wh}_1(X,x_0)$ is a topological group.
\end{proposition}
\proof
One only needs to show continuity of the concatenation operation. Let $\alpha$ and $\beta$ be two loops and $B([\alpha*\beta],U)$ be a neighborhood of their concatenation in $\pi^{wh}_1(X,x_0)$. The continuity of taking inverse (applied at the element $[\beta]^{-1}\in \pi^{wh}_1(X,x_0)$) implies the existence of a neighborhood $V$ of $x_0$ such that for any loop $\gamma$ in $V$ the inverse of $[\beta^{-1}*\gamma]$ belongs to $B([\beta],U)$, i.e. $[\beta^{-1}*\gamma]^{-1}=[\gamma^{-1}*\beta]=[\beta*\mu]$ for some loop $\mu$ in $U$.
Then for any $[\alpha*\gamma]\in B([\alpha],V)$ and any $[\beta*\delta]\in B([\beta],U)$ the concatenation $[\alpha*\gamma*\beta*\delta]=[\alpha*\beta*\mu*\delta]$ belongs to $B([\alpha*\beta],U)$.
\endproof

\begin{proposition}\label{DiscretnessOfBS}
$\pi_1^{wh}(X,x_0)$ is discrete if and only if $X$ is semi-locally simple-connected
at $x_0$.
\end{proposition}
\proof If there is a neighborhood $U$ of $x_0$ such that every loop
in $U$ at $x_0$ is null-homotopic in $X$, then $B([\alpha],U)$
contains only $[\alpha]$ among all homotopy classes of loops
at $x_0$. Conversely, if $B([\alpha],U)$ contains only $[\alpha]$ among all homotopy classes of loops
at $x_0$, then $[\alpha]=[\alpha\ast \gamma]$ for every loop $\gamma$ in $U$
at $x_0$, i.e. $[\gamma]=1$.

\begin{theorem}\label{FabelGeneralization}
Suppose $X$ is a path-connected space.
If $\pi_1^{top}(X,x_0)$ contains an isolated point, then $X$ is semi-locally simply-connected.
If $X$ is semi-locally simply-connected Peano space, then $\pi_1^{top}(X,x_0)$ is discrete.
\end{theorem}
\proof Choose an open set $V$ in $\widetilde X_{top}$ whose intersection
with $\pi_1^{top}(X,x_0)$ is exactly $[\alpha]$.
Given $x_1\in X$ choose a path $\lambda$ from $x_0$ to $x_1$.
Since $[(\alpha\ast\lambda)\ast \lambda^{-1}]\in V$, by \ref{BasicLemma}
there is a neighborhood $U$ of $x_1$ such that
$[(\alpha\ast\lambda)\ast \gamma\ast \lambda^{-1}]\in V$ for any
loop $\gamma$ in $U$ at $x_1$. Since $V$ contains only
the homotopy class of the loop $\alpha$,
$[\alpha]=[(\alpha\ast\lambda)\ast \gamma\ast \lambda^{-1}]$ and that implies
$\gamma$ is null-homotopic in $X$.

If $X$ is semi-locally simply-connected, then it is an SLT-space.
By~\ref{SLTImpliesPiContinuous}, we have $\pi_1^{top}(X,x)=\pi_1^{wh}(X,x)$ for every $x\in X$. Application of~\ref{DiscretnessOfBS} completes the proof.
\endproof

\begin{corollary}[Fabel \cite{Fab}]\label{FabelThm}
Suppose $X$ is a connected locally path-connected metrizable space.
The group $\pi_1^{top}(X,x_0)$ is discrete if and only if $X$ is semi-locally simply-connected.
\end{corollary}

\section{Uniform structures on $\widetilde X$}\label{SECTION-U-Structures}

In this section we introduce uniform structures on $\widetilde X$
and $\pi_1(X,x_0)$.

Given a uniform space $X$ and $x_0\in X$ we consider
the space $(X,x_0)^{(I,0)}$ of continuous functions with the {\bf uniform convergence structure} denoted by $(X,x_0)_{UC}^{(I,0)}$.
The base of that structure consists of elements $uc(E):=\{(\alpha,\beta) | (\alpha(t),\beta(t))\in E
\text{ for all }t\in I\}$, where $E$ is an entourage of $X$.

The projection
 $\pi\colon (X,x_0)^{(I,0)}\to \widetilde X$ is defined by $\pi(\alpha)=[\alpha]$.
 
Our first example of a uniform structure on $\widetilde X$
is the {\bf James uniform structure} (see \cite[p.120]{Jam}):

\begin{definition}\label{JamesStructureDef}
The {\bf James uniform structure} on $\widetilde X$ is generated by
sets $E^\ast$ consisting of pairs $([\alpha],[\beta])$ such that
$(\alpha,\beta)$ is a path in $E\subset X\times X$ originating at $(x_0,x_0)$.
\par $\widetilde X_J$ is the space $\widetilde X$ equipped with the James uniform structure.
\end{definition}

Since it is not entirely obvious $\{E^\ast\}$ form a base of a uniform structure,
let us provide some details that will be also useful later on.

\begin{definition}\label{PuncturedThingsDef}
By a {\bf punctured square} we mean a subset of $I\times I$
whose complement in $I\times I$ is a finite union of interiors
of disks $D_i$ such that $D_i\cap D_j=\emptyset$ if $i\ne j$ and $\bigcup D_i$
is contained in the interior of $I\times I$.
\par By an {\bf $E$-punctured-homotopy} from a path $\alpha$ to a path $\beta$
originating from the same point $x_0\in X$
we mean a map $H$ from 
a punctured square to $X$ satisfying the following conditions:
\begin{enumerate}
\item $H(t,0)=\alpha(t)$ and $H(t,1)=\beta(t)$ for all $t\in I$.
\item $H(0,s)=x_0$ for all $s\in I$.
\item $H(\{1\}\times I)$ is $E$-bounded.
\item $H(\partial D_i)$ is $E$-bounded for all disks $D_i$.
\end{enumerate}
\end{definition}

\begin{lemma}\label{PuncturedHomotopyExists}
Let $X$ be a uniform space with two given entourages $E$ and $F$
so that any pair $(x,y)\in F$ can be connected by an $E$-bounded path in $X$.
Suppose $\alpha,\beta\colon (I,0)\to (X,x_0)$ are two paths
in $X$. If the corresponding
elements $[\alpha],[\beta]$ of $\widetilde X$ satisfy $([\alpha],[\beta])\in \pi(uc(F))$, then
there is an $E^2$-punctured-homotopy from $\alpha$ to $\beta$.
\end{lemma}
\proof If $([\alpha],[\beta])\in \pi(uc(F))$, then there is $(\alpha^\prime,\beta^\prime)\in F$
so that $\alpha^\prime$ is homotopic to $\alpha$ rel. end-points
and $\beta^\prime$ is homotopic to $\beta$ rel. end-points.
Find $\epsilon >0$ with the property that $\alpha^\prime(J)$
and $\beta^\prime(J)$ are $E$-bounded for any subinterval $J$ of $I$
of diameter at most $\epsilon$. Subdivide $I$ into an increasing sequence
of points $t_0=0,t_1,\ldots,t_n=1$ so that $t_{i+1}-t_i < \epsilon$ for all $0\leq i < n$
and choose an $E$-bounded path $\gamma_i$ in $X$ from $\alpha^\prime(t_i)$
to $\beta^\prime(t_i)$. Use it to construct an $E^2$-punctured-homotopy
from $\alpha^\prime$ to $\beta^\prime$. Combine it with homotopies
from $\alpha$ to $\alpha^\prime$ and from $\beta$ to $\beta^\prime$
to create an $E^2$-punctured-homotopy from $\alpha$ to $\beta$.
\endproof

\begin{lemma}\label{IfPuncturedHomotopyExists}
Given two paths $\alpha,\beta\colon (I,0)\to (X,x_0)$, the corresponding
elements $[\alpha],[\beta]$ of $\widetilde X$ satisfy $([\alpha],[\beta])\in \pi(uc(E^2))$
if there is an $E$-punctured-homotopy $H$ from $\alpha$ to $\beta$.
\end{lemma}
\proof Without loss of generality we may assume the domain of $H$ is obtained from $I\times I$
by removing mutually disjoint disks $D_i$, $1\leq i\leq n$,
centered at $(x_i,\frac{1}{2})$ so that $x_i > x_j$ if $i < j$. 
Let $\alpha^\prime$ be the path in $X$ obtained from $H$
by starting at $(0,\frac{1}{2})$ and traversing upper semicircles of
boundries of disks $D_i$, reaching $(1,\frac{1}{2})$, and then going up
to $(1,1)$. Let $\beta^\prime$ be the path in $X$ obtained from $H$
by starting at $(0,\frac{1}{2})$ and traversing lower semicircles of
boundries of disks $D_i$, reaching $(1,\frac{1}{2})$, and then going down
to $(1,0)$. Notice $\alpha^\prime$ is homotopic to $\alpha$ rel. end-points,
$\beta^\prime$ is homotopic to $\beta$ rel. end-points,
and $(\alpha^\prime,\beta^\prime)\in uc(E^2)$.
\endproof

\begin{corollary}\label{GeneratingUCLemma}
If $X$ is uniformly locally
path-connected, then the projection
\par\noindent 
 $\pi\colon (X,x_0)^{(I,0)}\to \widetilde X$ $(\pi(\alpha)=[\alpha])$
generates the James uniform structure on $\widetilde X$.
\end{corollary}
\proof
First observe $E^\ast=\pi(uc(E))$, so if $\pi$ generates a uniform structure,
it must be identical with that of James.

Given an entourage $E$ of $X$ find an entourage $F$ such that $F^4\subset E$.
If $([\alpha],[\beta])\in \pi(uc(F))$ and $([\beta],[\gamma])\in \pi(uc(F))$
there exist $E^2$-punctured homotopies $H$ from $\alpha$ to $\beta$ and $G$
from $\beta$ to $\gamma$. The concatenation $H\ast G$ is an $E^2$-punctured-homotopy
from $\alpha$ to $\gamma$ hence $([\alpha],[\gamma])\in \pi(uc(F^4))\subset \pi(uc(E))$.
\endproof

Berestovskii and Plaut \cite{BP3} introduced a uniform structure
on the space of $E$-chains in $X$. It easily generalizes to a uniform structure on $\widetilde X$:

\begin{definition}\label{BPStructureDef}
The {\bf Berestovskii-Plaut uniform structure} on $\widetilde X$ (denoted by $\widetilde X_{BP}$) is generated by
sets $bp(E)$ consisting of pairs $([\alpha],[\beta])$ such that
$\alpha^{-1}\ast\beta$ is homotopic in $X$ rel. end-points to an $E$-bounded path.
\end{definition}
It is indeed a uniform structure: if $F^4\subset E$, then $(bp(F))^2\subset bp(E)$.
Note how much easier it is to introduce it in comparison to the
compact-open topology or James structure.

Notice that the
projection $\pi_X\colon \widetilde X_{BP}\to X$
is uniformly continuous. Also, it is surjective if and only if $X$ is path-connected.

\begin{proposition}\label{BPVsBogleySieradski}
If $X$ is completely regular, then the whisker
topology on $\widetilde X$ is identical with the topology induced on $\widetilde X$
by the Berestovskii-Plaut uniform structure corresponding to the canonical
uniform structure on $X$.
\end{proposition}
\proof Given a path $\alpha$ from $x_0$ to $x_1$ in $X$ and given
a neighborhood $U$ of $x_1$ in $X$, the set $B([\alpha],U)$
contains $B([\alpha],bp(E))$ for any entourage $E$ so that
$B(x_1,E^2)\subset U$. Also, $B([\alpha],bp(E))\supset B([\alpha],U)$
if $U\subset B(x_1,E)$.
\endproof

\begin{proposition}
The map $id\colon \widetilde X_{BP}\to \widetilde X_{J}$ is uniformly continuous for any uniform space $X$.
\end{proposition}
\proof 
If two paths $\alpha$ and $\beta$ have their homotopy classes close in the Berestovskii-Plaut uniform structure (i.e. $([\alpha],[\beta])\in bp(E)$ for some entourage $E$), then $\alpha^{-1}\ast\beta$ is homotopic in $X$ rel. end-points to some $E$-bounded path $\gamma$.
Thus $\beta$ is homotopic rel. end-points to $\alpha*\gamma$. Since $\gamma$ is $E$-bounded, the path $(\alpha*\gamma,\alpha*const)\subset X\times X$ is contained in $E$, therefore $[\alpha*\gamma]$ and $[\alpha*const]$ are $E^*$-close in the James uniform structure. Notice that $[\alpha*\gamma]=[\beta]$ and $[\alpha*const]=[\alpha]$.
\endproof

It is of interest to characterize spaces $X$ so that $id\colon \widetilde  X_{J}\to \widetilde X_{BP}$ is uniformly continuous.
Here is the corresponding concept to \ref{TopologicalSLTDef}
in the uniform category:

\begin{definition}\label{UniformSLTDef}
A uniform space $X$ is a {\bf uniform small loop transfer space} (a uniform SLT-space for short)
if for every entourage $E$ of $X$ there is an entourage $F$ of $X$
such that given two loops $\alpha\colon (S^1,1)\to (B(y,F),y)$
and $\beta\colon (S^1,1)\to (X,x)$ that are freely homotopic,
there is a loop $\gamma\colon (S^1,1)\to (B(x,E),x)$ that is homotopic
to $\beta$ rel. base point.
\end{definition}

\begin{theorem}\label{TwoStructuresOnTildeXComparison}
If $X$ is path-connected and uniformly locally path-connected, then
$id\colon \widetilde  X_{J}\to \widetilde X_{BP}$ is uniformly continuous
if and only if $X$ is a uniform SLT-space.
\end{theorem}
\proof Suppose $X$ is an SLT-space and $E$ is an entourage of $X$.
Find an entourage $F$ of $X$
such that given two loops $\alpha\colon (S^1,1)\to (B(y,F^2),y)$
and $\beta\colon (S^1,1)\to (X,x)$ that are freely homotopic,
there is a loop $\gamma\colon (S^1,1)\to (B(x,E),x)$ that is homotopic
to $\beta$ rel. base point.
It suffices to show $\pi(uc(F))\subset bp(E^3)$, so assume $([\alpha],[\beta])\in \pi(uc(F))$.
By \ref{PuncturedHomotopyExists} there is an $F^2$-punctured-homotopy $H$ from $\alpha$ to $\beta$.
Without loss of generality we may assume the domain of $H$ is obtained from $I\times I$
by removing mutually disjoint disks $D_i$, $1\leq i\leq n$,
centered at $(x_i,\frac{1}{2})$ so that $x_i > x_j$ if $i < j$. 
Let $\alpha_1$ be the path in $X$ obtained from re-parametrizing $H$
on the interval $[r_1,1]\times \{\frac{1}{2}\}$ traversed to the left (in the decreasing direction
of the first coordinate), where $(r_1,\frac{1}{2})$ belongs to $\partial D_1$.
More generally, $\alpha_i$, $n\ge i > 1$, is the path in $X$ obtained
from $H$ by traveling from the left-most point $L_{i-1}$ of $\partial D_{i-1}\cap I\times \{\frac{1}{2}\}$ to the right-most point $R_i$ of $\partial D_i\cap I\times \{\frac{1}{2}\}$.
Paths $\beta_i$, $1\leq i\leq n$, are obtained from $H$ by traveling
from $R_i$ to $L_i$ along the lower part of $\partial D_i$.
Paths $\gamma_i$, $1\leq i\leq n$, are obtained from $H$ by traveling
from $R_i$ to $L_i$ along the upper part of $\partial D_i$.
\par Define $l_k$ as $\prod\limits_{i=1}^k\alpha_i\ast \beta_i$
and define $u_k$ as $\prod\limits_{i=1}^k\alpha_i\ast \gamma_i$ for $k\leq n$.
We will show by induction on $k$ that $l_k\ast u^{-1}_k$ is homotopic
rel. end-points to a loop $\omega_k$ in $B(x_1,E)$, where $x_1=H(1,\frac{1}{2})$.
That will complete the proof of the first part of \ref{TwoStructuresOnTildeXComparison}
as $\alpha^{-1}\ast\beta$ is homotopic rel. end-points to 
$\delta\ast \omega_n\ast \mu$, where $\delta$ corresponds to $H|\{1\}\times [0,\frac{1}{2}]$
and $\mu$ corresponds to $H|\{1\}\times [\frac{1}{2},1]$.
\par Notice $l_1\ast u^{-1}_1$ is freely homotopic to $\beta_1\ast\gamma_1^{-1}$
that is contained in $B(R_1,F^2)$, so $\omega_1$ exists.
To do inductive step observe $l_k\ast u^{-1}_k\sim l_{k-1}\ast\beta_k\ast\gamma_k^{-1}\ast u_{k-1}^{-1}\sim (l_{k-1}\ast \beta_k\ast\gamma_k^{-1}\ast l_{k-1}^{-1})\ast
(l_{k-1}\ast u_{n-1}$ (here $\sim$ stands for homotopy rel. end-points).
Since $l_{k-1}\ast \beta_k\ast\gamma_k^{-1}\ast l_{k-1}^{-1}$ is freely
homotopic to a small loop $\beta_k\ast \gamma_k^{-1}$, it is homotopic
to a loop $\lambda$ contained in $B(x_1,E)$.
Thus $l_k\ast u^{-1}_k\sim \lambda\ast \omega_{k-1}$
and $\lambda\ast \omega_{k-1}$ is contained in $B(x_1,E)$.
\par Suppose $id\colon \widetilde  X_{UC}\to \widetilde X_J$ is uniformly continuous
and $E$ is an entourage of $X$. Pick an entourage $F$ of $X$
such that $uc(F^2)\subset bp(E)$.
Given a loop $\alpha$ in $B(x_1,F)$ and a path $\lambda$ from $x_1$
to $x_0$ represent $\alpha$ as $\beta\ast\gamma$
and notice $([\beta^{-1}\ast \lambda],[\gamma\ast \lambda])\in uc(F^2)$.
Hence $([\beta^{-1}\ast \lambda],[\gamma\ast \lambda])\in bp(E)$
and there is a loop $\mu$ in $B(x_0,E)$ such that
$\mu$ is homotopic rel. end-points to 
$(\beta^{-1}\ast \lambda)^{-1}\ast(\gamma\ast \lambda)\sim \lambda^{-1}\ast\alpha\ast\lambda$. That proves any loop freely homotopic to $\alpha$
is homotopic rel. end-points to a small loop.
\endproof

\begin{problem}\label{AreSLTsPoincare}
Is there a uniform SLT-space $X$ that is not a uniform Poincare space?
\end{problem}

We introduce the lasso uniformity on the set $\widetilde X$ for any uniform space $X$. 

\begin{definition}
Let $E$ be an entourage of a uniform space $X$ and $x$ be a point in $X$. A path $l$ is called {\it $E$-lasso based at the point $x$} if $l$ is equal to a finite concatenation of loops $\alpha_n\ast\gamma_n\ast\alpha_n^{-1}$,
where each loop $\gamma_n$ is $E$-bounded
and $\alpha_n$ is a path from $x$ to $\gamma_n(0)$.
\end{definition}

\begin{definition}
The {\bf lasso uniform structure} on $\widetilde X$ (denoted by $\widetilde X_{l}$) is generated by
sets $l(E)$ consisting of pairs $([\alpha],[\beta])$ such that for some $E$-lasso $l$ the path
$\alpha^{-1}*l\ast\beta$ is homotopic in $X$ rel. end-points to an $E$-bounded path.
\end{definition}

\begin{lemma}\label{lassoClosedMeansPuncturedHomotopic}
Let $X$ be a uniform space, $E$ be an entourage, and $\alpha$ and $\beta$ be two paths in $X$.
The classes $[\alpha]$ and $[\beta]$ are $E$-close in the lasso uniformity on $\widetilde X$ if and only if there is an $E$-punctured homotopy from $\alpha$ to $\beta$.
\end{lemma}
\proof 
We may assume that the domain of an $E$-punctured homotopy $H$ from $\alpha$ to $\beta$ is obtained from $I\times I$
by removing mutually disjoint congruent disks $D_i$, $1\leq i\leq n$,
centered at $(x_i,\frac{1}{2})$ so that $x_i > x_j$ if $i < j$. 
Let $\alpha_i$ be the straight path from the point $(0,0)$ to the point $p_i$ of the boundary of the disk $D_i$ with the smallest $y$-coordinate (so that any path $\alpha_i$ is disjoint from other paths $\alpha_j$ and disks $D_j$). Let the loop $\gamma_i$ be based at the point $p_i$ and traversing the boundary of the disk $D_i$ once. Then $H$ makes the path $\prod\limits_{i=1}^n\alpha_i\ast \gamma_i*\alpha_i^{-1}$ an $E$-lasso $l$ in $X$ and thus $H$ can be interpreted as a homotopy between $\alpha^{-1}*l\ast\beta$ and an $E$-bounded path. Conversely, a homotopy between $\alpha^{-1}*l\ast\beta$ and an $E$-bounded path can be interpreted as an $E$-punctured homotopy from $\alpha$ to $\beta$.
\endproof

The universality condition of the universal Peano space $P(X)$ (see~\ref{ULPCspacesDef}) allows us to identify the sets $\widetilde X$ and $\widetilde {P(X)}$ since any path in $P(X)$ projects to a path in $X$ and any path in $X$ lifts uniquely to a path in $P(X)$.

\begin{proposition}\label{X_L=PX_J}
Let $X$ be a uniform space. The natural identification of the sets $\widetilde X$ and $\widetilde {P(X)}$ generates a uniform homeomorphism of $\widetilde X_l$ and $\widetilde {P(X)}_J$
\end{proposition}
\proof 
Suppose that two paths $\alpha$ and $\beta$ in $X$ are $E$-close in the lasso uniformity. By~\ref{lassoClosedMeansPuncturedHomotopic}, there is an $E$-punctured homotopy from $\alpha$ to $\beta$.
By~\ref{IfPuncturedHomotopyExists}, the classes $[\alpha]$ and $[\beta]$ have uniformly $E^2$-close representatives (notice that the representatives constructed in the proof of~\ref{IfPuncturedHomotopyExists} are also uniformly $E^2$-close in $P(X)$). By~\ref{GeneratingUCLemma}, the classes $[\alpha]$ and $[\beta]$ are close in the James uniform structure on $\widetilde{P(X)}$.

Suppose that two paths $\alpha$ and $\beta$ in $P(X)$ are $pc(E)$-close in the James uniform structure on $\widetilde{P(X)}$. By~\ref{PuncturedHomotopyExists}, there is an $E^2$-punctured-homotopy from $\alpha$ to $\beta$.
By~\ref{lassoClosedMeansPuncturedHomotopic}, the classes $[\alpha]$ and $[\beta]$ are $E^2$-close in the lasso uniformity on $\widetilde X$.
\endproof

\begin{corollary}\label{X_l=X_J}
For any connected and uniformly locally path-connected uniform space $X$ the James uniformity on $\widetilde X$ is identical with the lasso uniformity.
\end{corollary}

\proof 
For any connected and uniformly locally path-connected uniform space $X$ we have $P(X)=X$.
\endproof

A uniform structure on $\widetilde X$ induces a uniform structure on the fundamental group $\pi_1(X,x_0)$. It is natural to ask if this uniform structure makes the group a topological group. 

Recall that the quotient of the compact-open topology does not make $\pi_1(X,x_0)$ a topological group for the single reason of the product of quotient maps not being a quotient map~\cite{Fab3}. If the space $X$ is a uniform space, the situation is better:  each product of uniform quotient mappings is a uniform quotient mapping~\cite{HR}. Thus the fundamental group equipped with the James uniform structure is a topological group. The Berestovskii-Plaut uniform structure is unlikely to make the fundamental group a topological group for the same reason as the whisker topology does not make $\pi_1(X,x_0)$ a topological group (recall~\ref{BPVsBogleySieradski}). 

Although the lasso uniformity on $\pi_1(X,x_0)$ can be related to the James uniformity, we give a direct simple proof of the following:

\begin{proposition}
The lasso uniform structure makes the fundamental group a topological group.
\end{proposition}
\proof 
By definition, two elements of $[\alpha], [\beta]\in\pi_1(X,x_0)$ are $E$-close if for some $E$-lasso $l$ the loop
$\alpha^{-1}*l\ast\beta$ is homotopic to an $E$-bounded loop $\gamma$. That means for the $E$-lasso $L=l*\beta*\gamma*\beta^{-1}$ the loop $\alpha^{-1}*L\ast\beta$ is homotopically trivial. Thus two elements of $[\alpha], [\beta]\in\pi_1(X,x_0)$ are $E$-close if the loop $\alpha*\beta^{-1}$ is homotopic to an $E$-lasso.

The inverse operation in $\pi_1(X,x_0)$ is continuous in the lasso uniformity because if the loop $\alpha*\beta^{-1}$ is homotopic to an $E$-lasso $L$, then the loop $(\alpha^{-1})*(\beta^{-1})^{-1}=\alpha^{-1}*\beta$ is homotopic to an $E$-lasso $\alpha^{-1}*L^{-1}*\alpha$ (use~\ref{ConjugateLassoIsLasso}).

The product operation in $\pi_1(X,x_0)$ is continuous in the lasso uniformity because if the loop $\alpha*\beta^{-1}$ is homotopic to an $E$-lasso $L$ and the loop $\gamma*\delta^{-1}$ is homotopic to an $E$-lasso $M$, then the loop $(\alpha*\gamma)*(\beta*\delta)^{-1}\sim\alpha*M*\beta^{-1} \sim\alpha*\beta^{-1}*\beta*M*\beta^{-1} \sim L*(\beta*M*\beta^{-1})$ is homotopic to an $E$-lasso (use~\ref{ConjugateLassoIsLasso}).
\endproof


\begin{thebibliography}{99}

\bibitem{AT}
A. ArhangelÕskii, M. Tkachenko, {\em Topological Groups and Related Structures}, Atlantis Studies in Mathematics, 1. Atlantis Press, Paris, 2008.

\bibitem{BP3}
V. Berestovskii, C. Plaut, {\em Uniform universal covers of uniform spaces},
Topology Appl. 154 (2007), 1748--1777.


\bibitem{BogSie}  W.A. Bogley, A.J. Sieradski, {\em Universal path spaces}, http://oregonstate.edu/\symbol{126}bogleyw/\#research

\bibitem{Bra} J. Brazas,
{\em The fundamental group as a topological group}, Preprint arXiv:1009.3972v5.

\bibitem{BDLM1} N. Brodskiy, J. Dydak, B. Labuz, A. Mitra,
{\em Rips complexes and covers in the uniform category}, Preprint math.MG/0706.3937, accepted by the Houston Journal of Math.

\bibitem{BDLM2} N. Brodskiy, J. Dydak, B. Labuz, A. Mitra,
{\em Group actions and covering maps in the uniform category}, Topology Appl., 157 (2010), 2593--2603. doi:10.1016/j.topol.2010.07.011.

\bibitem{BDLM3} N. Brodskiy, J. Dydak, B. Labuz, A. Mitra,
{\em Covering maps for locally path-connected spaces}, Preprint arXiv:0801.4967.

\bibitem{Dug} J. Dugundji, {\em A topologized fundamental group}, 
Proc. Nat. Acad. Sci. U. S. A. 36, (1950). 141--143.

\bibitem{Dyd} J. Dydak, {\em Partitions of unity}, Topology Proceedings
  27 (2003),  125--171.
  
  \bibitem{DydSeg}
J. Dydak, J. Segal, {\em Shape theory: An introduction},
Lecture Notes in Math. 688, 1--150, Springer Verlag 1978.

\bibitem{Fab1}  P. Fabel, {\em The Hawaiian earring group and metrizability},
ArXiv:math.GT/0603252.

\bibitem{Fab} P. Fabel, {\em Metric spaces with discrete topological fundamental group}, Topology and its Applications 154 (2007), 635--638.

\bibitem{Fab3}  P.~Fabel, {\em Multiplication is Discontinuous in the Hawaiian Earring Group (with the Quotient Topology)},   Bull. Pol. Acad. Sci. Math. 59 (2011), no. 1, 77Ð83. DOI: 10.4064/ba59-1-09

\bibitem{FisZas}  H. Fischer, A. Zastrow, {\em Generalized universal
coverings and the shape group}, Fundamenta Mathematicae 197 (2007), 167--196.

\bibitem{GhyHar}
E. Ghys and E. de la Harpe (eds), {\em Sur les groupes hyperboliques d'apres de Mikhael
Gromov}, Birkh\" auser, 1990.

\bibitem{Hur} W. Hurewicz, {\em Homotopie, Homologie und lokaler Zusammenhang}, Fundamenta Mathematicae 25 (1935), 467--485.

\bibitem{HR} M. Hu\v{s}ek, M. Rice, {\em  Productivity of coreflective subcategories of uniform spaces},  General Topology Appl. 9 (1978), no. 3, 295--306.
  
\bibitem{Isb}
J. Isbell, {\em Uniform Spaces}, Mathematical Surveys, vol. 12, American Mathematical Society, Providence, RI, 1964.

\bibitem{Jam}
I.M. James, {\em Introduction to Uniform Spaces}, London
Math. Soc. Lecture Notes Series 144, Cambridge
University Press, 1990.
 
 \bibitem{Kel}  J.L. Kelley,
{\em General topology},
D. Van Nostrand Company, Inc., Toronto-New York-London, 1955.

\bibitem{KraMin}
J. Krasinkiewicz, P. Minc,
{\em Generalized paths and pointed 1-movability}, Fund. Math. 104 (1979), no. 2, 141--153.

\bibitem{Lab}
B. Labuz, {\em Inverse limits and uniform covers},
in preparation.

\bibitem{Lim}
E.L. Lima, {\em Fundamental groups and covering spaces},
AK Peters, Natick, Massachusetts, 2003.

\bibitem{MarSeg}
S. Marde\v si\' c, J. Segal,
 {\em Shape theory}, North-Holland Publ.Co., Amsterdam 1982.

\bibitem{Mel} S.A. Melikhov, {\em Steenrod homotopy}, (Russian) Uspekhi Mat. Nauk 64 (2009), no. 3(387), 73--166; translation in Russian Math. Surveys 64 (2009), no. 3, 469--551.

\bibitem{Mun}
J.R. Munkres, {\em Topology}, Prentice Hall, Upper Saddle River, NJ 2000.

\bibitem{Nakano}
H. Nakano, {\em Uniform spaces and transformation groups}. Wayne State University Press, Detroit, Mich. 1968 xv+253 pp.

\bibitem{Paw} J. Pawlikowski, {\em The fundamental group
of a compact metric space},
Proceedings of the
American Mathematical Society, 126 (1998), 3083--3087.

\bibitem{Pla}
C. Plaut, {\em Quotients of uniform spaces},
Topology Appl. 153 (2006), 2430--2444.

\bibitem{Spa}
E. Spanier,
{\em Algebraic topology}, McGraw-Hill, New York 1966.

\bibitem{V1}
\v{Z}. Virk, {\em Small loop spaces},
Topology Appl. 157 (2010), 451--455.

\bibitem{V2}
\v{Z}. Virk, {\em Homotopical smallness and closeness},
Topology Appl. 158 (2011), 360--378.

\end{thebibliography}
\end{document}